\documentclass[12pt,leqno]{article}
\usepackage{amsfonts}
\usepackage{amsmath, amsthm, amsfonts, amssymb, color}
\usepackage{mathrsfs}
\usepackage{color}
\theoremstyle{definition}

\newcommand{\scr}[1]{\mathscr #1}
\definecolor{wco}{rgb}{0.5,0.2,0.3}

\numberwithin{equation}{section} \theoremstyle{remark}

\newcommand{\ua}{\uparrow}

\title{{\bf  Integration by Parts Formula and Applications  for SDEs with L\'evy Noise}\footnote{Supported in
 part by  NNSFC(11131003), SRFDP, the Fundamental Research Funds for the Central Universities.} }
\author{
{\bf     Feng-Yu Wang  }\\
\footnotesize{ School of Mathematical Sciences,
Beijing Normal
University, Beijing 100875, China}\\
 \footnotesize{ Department of Mathematics,
Swansea University, Singleton Park, SA2 8PP, United Kingdom}\\
\footnotesize{  wangfy@bnu.edu.cn, F.-Y.Wang@swansea.ac.uk}}
\begin{document}
\allowdisplaybreaks
\def\R{\mathbb R}  \def\ff{\frac} \def\ss{\sqrt} \def\B{\mathbf
B}
\def\N{\mathbb N} \def\kk{\kappa} \def\m{{\bf m}}
\def\ee{\varepsilon}\def\ddd{D^*}
\def\dd{\delta} \def\DD{\Delta} \def\vv{\varepsilon} \def\rr{\rho}
\def\<{\langle} \def\>{\rangle} \def\GG{\Gamma} \def\gg{\gamma}
  \def\nn{\nabla} \def\pp{\partial} \def\E{\mathbb E}
\def\d{\text{\rm{d}}} \def\bb{\beta} \def\aa{\alpha} \def\D{\scr D}
  \def\si{\sigma} \def\ess{\text{\rm{ess}}}
\def\beg{\begin} \def\beq{\begin{equation}}  \def\F{\scr F}
\def\Ric{\text{\rm{Ric}}} \def\Hess{\text{\rm{Hess}}}
\def\e{\text{\rm{e}}} \def\ua{\underline a} \def\OO{\Omega}  \def\oo{\omega}
 \def\tt{\tilde} \def\Ric{\text{\rm{Ric}}}
\def\cut{\text{\rm{cut}}} \def\P{\mathbb P} \def\ifn{I_n(f^{\bigotimes n})}
\def\C{\scr C}      \def\aaa{\mathbf{r}}     \def\r{r}
\def\gap{\text{\rm{gap}}} \def\prr{\pi_{{\bf m},\varrho}}  \def\r{\mathbf r}
\def\Z{\mathbb Z} \def\vrr{\varrho} \def\ll{\lambda}
\def\L{\scr L}\def\Tt{\tt} \def\TT{\tt}\def\II{\mathbb I}
\def\i{{\rm in}}\def\Sect{{\rm Sect}}  \def\H{\mathbb H}
\def\M{\scr M}\def\Q{\mathbb Q} \def\texto{\text{o}} \def\LL{\Lambda}
\def\Rank{{\rm Rank}} \def\B{\scr B} \def\i{{\rm i}} \def\HR{\hat{\R}^d}
\def\to{\rightarrow}\def\l{\ell}\def\iint{\int}
\def\EE{\scr E}
\def\A{\scr A}
\def\BB{\scr B}\def\Ent{{\rm Ent}}

\maketitle

\begin{abstract}  By using the Malliavin calculus and finite-jump approximations, the Driver-type integration by parts formula is established for the semigroup associated to stochastic differential equations with  noises containing a subordinate Brownian motion. As applications, the shift-Harnack inequality and heat kernel estimates are derived. The main results are illustrated by SDEs driven by $\aa$-stable like processes.
\end{abstract} \noindent
 \vskip 2cm

\section{Introduction}

A significant application of the Malliavin calculus is to describe the density of a Wiener functional using the integration by parts formula.
In 1997, Driver \cite{Driver} established the following integration by parts formula for the heat semigroup $P_t$ on a compact Riemannian manifold $M$:
$$P_t (\nn_Z f)= \E\{f(X_t)N_t\},\ \ f\in C^1(M), Z\in \scr X,$$ where $\scr X$ is the set of all smooth vector fields on $M$, and $N_t$ is a random variable depending on $Z$ and the curvature tensor. From this formula we are able to characterize   the derivative w.r.t. the second variable $y$ of the heat kernel $p_t(x,y)$, see \cite{W12} for a recent study on integration by parts formulas and applications for stochastic differential equations driven by Wiener processes. The backward coupling method developed in \cite{W12} has been also used in \cite{Fan,Zhang} for SDEs driven by fractional Brownian motions and SPDEs driven by Wiener processes.  The purpose of this paper is to investigate the integration by parts formula and applications for SDEs  driven by purely jump L\'evy noises.

Consider the following stochastic equation on $\R^d$:
\beq\label{1.1} X_t= X_0+ \int_0^t b_s(X_s)\d s+\int_0^t \si_s \d W_{S(s)} +V_t,\ \ t\ge 0,\end{equation}
where $$\si: [0,\infty)\to\R^d\otimes\R^d,\ \ b: [0,\infty)\times \R^d\to \R^d$$ are measurable and  locally bounded,
 $W:=(W_t)_{t\ge 0}, S:=(S(t))_{t\ge 0}$ and $V:=(V_t)_{t\ge 0}$ are independent stochastic processes such that

\beg{enumerate} \item [(i)] $W$ is the Brownian motion on $\R^d$ with $W_0=0$;
\item[(ii)]  $V$ is a c\'adl\'ag
  process on $\R^d$ with $V_0=0$;
  \item[(iii)] $S$ is the subordinator induced by a Bernstein function $B$, i.e. $S$ is a  one-dimensional   increasing L\'evy process with $S(0)=0$ and Laplace transform
$$\E \e^{-r S(t)}= \e^{-t B(r)},\,\,\ \ t,r\ge 0.$$ \end{enumerate}
Then $(W_{S(t)})_{t\ge 0}$ is a the L\'evy process known as the subordinate Brownian motion with subordinator $S$,  see e.g.  \cite{A,J}.

For this equation the Bismut formula and Harnack inequalities  have been studied in \cite{Z} and \cite{WW} by using regularization approximations of $S(t)$, but the study of the integration by parts formula and shift-Harnack inequality is not yet done.

To establish the integration by parts formula, we need the following assumptions.
\beg{enumerate} \item[(H1)]   $b_t\in C^2(\R^d)$  such that for some increasing $K_1,K_2\in C([0,\infty))$,
$$\|\nn b_t\|_\infty:=\sup_{x\in\R^d} \|\nn b_t(x)\|\le K_1(t),\ \ \|\nn^2 b_t\|_\infty:=\sup_{x\in\R^d}\| \nn^2 b_t(x)\|\le K_2(t),\ \ t\ge 0.$$
\item[(H2)] $\si_t$ is invertible   such that for some increasing $\ll_1,\ll_2\in C([0,\infty)),$
$$\|\si_t\|\le \ll_1(t),\ \ \|\si_t^{-1}\|\le \ll_2(t),\ \ t\ge 0.$$
\end{enumerate}

 By (H1) and (H2), for any   $x\in \R^d$  the equation (\ref{1.1}) with $X_0=x$ has a unique solution $X_t(x)$. Let $P_t$ be the associated Markov semigroup, i.e.
$$P_t f(x)= \E f(X_t(x)),\ \ f\in\B_b(\R^d), t\ge 0, x\in\R^d.$$

As already observed in \cite{W12} that comparing with the Bismut formula,   the integration by parts formula is usually harder to establish.
 To strengthen this observation, we explain below   that the regularization argument used in \cite{Z} for the Bismut formula is no longer valid for the integration by parts formula.
For simplicity, let us consider the case that $V_t=0, b_t=b$ and $\si_t=\si$.  As in \cite{Z}, for any $\vv>0$ let
$$S_\vv(t)= \ff 1 \vv\int_t^{t+\vv} S(s)\d s+\vv t,\ \ t\ge 0.$$ Then $S_\vv(\cdot)$ is differentiable and $S_\vv\downarrow S$ as $\vv\downarrow 0.$ Consider the equation (note that we have assumed $V_t=0, b_t=b$ and $\si_t=\si$)
$$\d X_t^\vv= b(X_t^\vv)\d t +\si \d W_{S_\vv(t)},\ \ X_0^\vv = X_0.$$
To make use the existing derivative formulas for SDEs driven by the Brownian motion, we take
$Y_t^\vv= X_{S_\vv^{-1}(t)}^\vv$ so that this equation reduces to
$$\d Y_t^\vv= b (Y_t^\vv) (S_\vv^{-1})'(t) \d t +\si \d W_t,\ \ Y_0^\vv=X_0.$$ In \cite{Z}, by using a known Bismut formula for $Y_t^\vv$ and letting $\vv\to 0$, the corresponding formula for $X_t$ is established. The crucial point for this argument is that the Bismut formula for $Y_t^\vv$   converges as $\vv\to 0$. However, since $S$ is not differentiable,   the existing integration by parts formula of $Y_t^\vv$ (see e.g. \cite[Theorem 5.1]{W12} with $H=\R^d$ and $ A=0$)
 $$\E(\nn_v f)(Y_T^\vv) =\ff 1 T \E\bigg\{f(Y_T^\vv)\int_0^T\<\si^{-1}(v-t(S_\vv^{-1})'(t)\nn_v b(Y_t^\vv), \d W_t\>\bigg\}$$
   does not converge to any explicit formula as $\vv\to 0$, except $\nn_v b$ is trivial.

So, to establish the integration by parts formula, we will take a different approximation argument, i.e. the finite-jump approximation   used in \cite{WXZ} to establish the Bismut formula for   SDEs with  multiplicative  L\'evy noises. We have to indicate that in this paper we are not able to establish the integration by parts formula for SDEs with multiplicative  L\'evy noises. Note that even for SDEs driven by multiplicative Gaussian noises, the existing integration by parts formula using the Malliavin covariant matrix is in general less explicit.

 To state our main result, we introduce the $\R^d\otimes \R^d$-valued process $J_t$, which solves the ODE
 \beq\label{*0} \ff{\d}{\d t} J_t= (\nn b_t)(X_t) J_t,\ \ J_0=I,\end{equation} where for $t\ge 0$ and $x\in\R^d, (\nn b_t)(x)\in \R^d\otimes\R^d$ is determined by
 $$(\nn b_t)(x)v = (\nn_v b_t)(x),\ \ v\in \R^d.$$
 By (H1),   both $J_t$ and $J_t^{-1}$ are locally bounded in $t$.

\beg{thm}\label{T1.1}  Assume {\rm (H1)} and {\rm (H2)}. Let $T>0.$  If
\beq\label{S} \E S(T)^{-\ff 12} <\infty,\end{equation} then for any $v\in \R^d$ and $f\in C_b^1(\R^d)$,
\beq\label{I} P_T(\nn_v f) = \E\Big\{f(X_T) \ff{M_T^v}{S(T)}\Big\},\end{equation} where
\beg{equation*}\beg{split} M_T^v:=  &\bigg\<\int_0^T  \big(\si_t^{-1} J_t\big)^*\d W_{S(t)},  J_T^{-1}v\bigg\>\\
&+\int_0^T\d S(t)\int_t^T {\rm Tr}\Big\{\si_t^{-1}J_tJ_s^{-1}\big(\nn\nn_{J_sJ_T^{-1}v}b_s\big)(X_s) J_sJ_t^{-1}\si_t\Big\}\d s\end{split}\end{equation*}  is in $L^1(S(T)^{-1}\d \P)$. Consequently, $P_T$ has density $p_T(x,y)$ with respect to the Lebesgue measure, which is differentiable in $y$  with
 $$\nn_v\log p_T(x,\cdot)(y)= -\E\Big(\ff{M_T^v}{S(T)}\Big|X_T(x)=y\Big),\ \ v,x\in\R^d.$$   \end{thm}

Below we present some consequences of Theorem \ref{T1.1} concerning derivative estimates and shift-Harnack inequalities.  For non-negative $f\in \B_b(\R^d)$ and $T>0$, let
$${\rm Ent}_{P_T}(f)= P_T(f\log f)- (P_Tf)\log P_Tf$$ be the relative entropy of $f$ with respect to $P_T$.

\beg{cor} \label{C1.2} Assume that {\rm (H1)}, {\rm (H2)} and     $\eqref{S}$ hold.
Let $$\bb(T)=  d T\ll_1(T) \ll_2(T)K_2(T)\e^{3TK_1(T)},\ \ T>0.$$
\beg{enumerate} \item[$(1)$] For any $T>0$ and $ v\in\R^d$,
\beg{equation*}\beg{split} &\| P_T(\nn_v f)\|_\infty\le |v|\cdot \|f\|_\infty \Big(\ll_2(T)\e^{TK_1(T)} \E S(T)^{-\ff 1 2} + \bb(T)\Big),\ \  f\in C_b^1(\R^d),\\
&  \int_{\R^d} |\nn_v p_T(x,\cdot)|(y) \d y\le |v| \Big(\ll_2(T)\e^{TK_1(T)} \E S(T)^{-\ff 1 2} + \bb(T)\Big),\  \ x\in\R^d.\end{split}\end{equation*}
\item[$(2)$] For any $p>1$ there exists a constant $C(p)\ge 1$ such that for any $T>0$,
\beg{equation*}\beg{split} &| P_T(\nn f)|\le C(p) (P_T|f|^p)^{\ff 1 p} \Big(\ll_2(T)\e^{TK_1(T)} \big(\E S(T)^{-\ff p {2(p-1)}}\big)^{\ff {p-1}p}
  + \bb(T)\Big),\ f\in C_b^1(\R^d),\\
 &  \int_{\R^d} |\nn\log p_T(x,\cdot)|^{\ff p{p-1}}(y)p_T(x,y)\d y \\
  &\qquad  \le C(p) \Big(\ll_2(T)\e^{TK_1(T)} \big(\E S(T)^{-\ff p {2(p-1)}}\big)^{\ff {p-1}p}
  + \bb(T)\Big),\ \ x\in\R^d. \end{split}\end{equation*}
\item[$(3)$] For any $\dd>0,v,x\in\R^d$ and positive $f\in \B_b(\R^d)$,
\beg{equation*}\beg{split} & |P_T(\nn_v f)|\le  \dd{\rm Ent}_{P_T}(f)
+   (P_Tf) \Big(\bb(T)|v|
 + \dd\log\E\exp\Big[\ff{\ll_2(T)^2|v|^2\e^{2TK_1(T)}}{2\dd^2 S(T)}\Big]\Big),\\
 &  \int_{\R^d} \exp\Big[\ff{|\nn_v\log p_T(x,\cdot)|(y)}\dd\Big]p_T(x,y)\d y
 \le \E\exp\Big[\ff{\bb(T)|v|}\dd +\ff{\ll_2(T)^2|v|^2\e^{2TK_1(T)}}{2\dd^2S(T)} \Big].\end{split}\end{equation*} \end{enumerate}
\end{cor}

\beg{cor} \label{C1.3} Assume {\rm (H1)} and {\rm (H2)}. Let $p>1, T>0.$ If
$$\GG_{T,p}(r):= \E\exp\bigg[\ff{p^2\ll_2(T)^2\e^{2TK_1(T)}r^2}{2(p-1)^2S(T)}\bigg]<\infty,\ \   r\ge 0,$$ then the shift-Harnack inequality
\beq\label{SH}(P_Tf)^p(x)\le \exp\Big[ \ff{p(\log p)\bb(T)|v|}{p-1}
+\ff{p-1}p\log \GG_{T,p}(|v|)\Big]P_T(f^p(v+\cdot))(x)\end{equation} holds for all $v,x\in\R^d$ and positive $f\in \B_b(\R^d).$ Consequently,
$$\sup_{x\in\R^d} \int_{\R^d} p_T(x,y)^{\ff p{p-1}}\d y \le \bigg(\int_{\R^d}\exp\Big[ -\ff{p(\log p)\bb(T)|v|}{p-1}
-\ff{p-1}p\log \GG_{T,p}(|v|)\Big]\d v\bigg)^{\ff{-1}{p-1}}.$$ \end{cor}

\

To illustrate the above results,   we consider below the SDE driven by $\aa$-stable like noises.

\

\beg{cor}\label{C1.4} Assume {\rm (H1)} and {\rm (H2)}. Let $B(r)\ge cr^{\ff\aa 2}$ for $r\ge r_0$, where  $\aa\in (0,2)$ and  $c,r_0>0$ are constants.
\beg{enumerate} \item[$(1)$] For any $p>1$ there exists a constant $C(p)>0$ such that
\beg{equation*}\beg{split} &|P_T(\nn f)|\le \ff{C(p)(P_T|f|^p)^{\ff 1 p}}{1\land T^{\ff 1 \aa}},\ \ T>0, f\in C_b^1(\R^d),\\
 &\sup_{x\in\R^d} \int_{\R^d} |\nn\log p_T(x,\cdot)|^{\ff p{p-1}}(y)p_T(x,y)\d y\le \ff{C(p)}{1\land T^{\ff 1 \aa}},\ \ T>0.\end{split}\end{equation*} \item[$(2)$] Let $\aa\in (1,2).$ Then there exists a constant $C>0$ such that for any $p>1,\dd>0, v\in\R^d$ and $f\in C^1(\R^d)$,
 \beg{equation*}\beg{split}&|P_T(\nn_v f)|\le \dd{\rm Ent}_{P_T}(f) + (P_Tf) \Big(\bb(T)|v|+\ff{C|v|^2}{\dd^2(1\land T)^{\ff 2\aa}}+ \ff{C|v|^{\ff\aa{\aa-1}} }{\{\dd^\aa (1\land T)\}^{\ff 1 {\aa-1}}}\Big),\\
 &\sup_{x\in\R^d}\int_{\R^d} \exp\Big[\ff{|\nn_v\log p_T(x,\cdot)(y)|}{\dd}\Big]p_T(x,y)\d y\\
 &\qquad \le \exp\Big[\bb(T)|v|+ \ff{C|v|^2}{\dd^2(1\land T)^{\ff 2\aa}}+ \ff{C|v|^{\ff\aa{\aa-1}} }{\{\dd^\aa (1\land T)\}^{\ff 1 {\aa-1}}}\Big].\end{split}\end{equation*}
 \item[$(3)$] Let $\aa\in (1,2).$ Then there exists a constant $C>0$ such that for any $p>1,T>0, v\in\R^d$ and positive $f\in\B_b(\R^d), $
 \beg{equation*}\beg{split} &(P_Tf)^p\le \exp\Big[\ff{C(p\log p)|v|}{p-1}+\ff{Cp|v|^2}{(p-1)(1\land T)^{\ff 2\aa}}+ \ff{Cp^{\ff 1{\aa-1}}|v|^{\ff\aa{\aa-1}}}{[(p-1)(1\land T)]^{\ff 1{\aa-1}}}\Big]P_T( f^p(v+\cdot)),\\
 &\sup_{x\in\R^d} \int_{\R^d} p_T(x,y)^{\ff p{p-1}}\d y\le \ff{1}{(1\land T)^{\ff d{\aa(p-1)}}}\exp\Big[\ff{Cp\log p}{(p-1)^2}+\ff{Cp^{\ff 1 {\aa-1}}}{(p-1)^{\ff\aa{\aa-1}}}\Big].\end{split}\end{equation*} \end{enumerate}\end{cor}

 In the next section we will prove the integration by parts formula for time-changes of finite jumps, from which we will be able to present in Section 3  complete proofs of the above results.

 \section{Integration by parts formula with finite-jump}

 In this section, we let $\ell$ be a c\'adl\'ag and increasing function on $[0,\infty)$ with $\ell(0)=0$ such that the set $\{t\in [0,T]: \DD \ell(t):= \ell(t)-\ell(t-)>0\}$ is finite. We call $\ell $ a path of $S$ with finite many jumps on $[0,T]$. Let $X_t^\ell$ solve the equation
\beq\label{2.10}   X_t^\ell =X_0 + \int_0^t b_s(X_s^\ell)\d s+ \int_0^t\si_s\d W_{\ell(s)}+V_t,\ \ t\ge 0,\end{equation}   and let $P_t^\ell$ be the associated Markov semigroup; i.e.
 $$P_t^\ell f:= \E f(X_t^\ell).$$ Moreover, let $J_t^\ell$ solve the ODE
\beq\label{*1.2}\ff{\d}{\d t} J_t^\ell= (\nn b_t)(X_t^\ell) J_t^\ell,\ \ J_0^\ell=I.\end{equation}
 The main result in this section is the following.

 \beg{thm}\label{T2.1} Let $\ell$ be a path of $S$ with finite many jumps on $[0,T].$ Then
 $$P_T^\ell (\nn_vf)= \E\big\{f(X_T^\ell) M_T^{\ell,v}\big\},\ \ v\in\R^d, f\in C_b^1(\R^d),$$ where
 \beg{equation*}\beg{split} M_T^{\ell, v}:= &\bigg\<\int_0^T \big(\si_t^{-1} J_t^\ell\big)^*\d W_{\ell(t)}, (J_T^\ell)^{-1}v \bigg\>\\
 &+ \int_0^T \d \ell(t) \int_t^T {\rm Tr}\Big\{ \si_t^{-1} J_t^\ell(J_s^\ell)^{-1}\big(\nn\nn_{J_s^\ell(J_t^\ell)^{-1}v} b_s\big)(X_s^\ell)J_s^\ell(J_t^\ell)^{-1}\si_t \Big\}\d s.\end{split}\end{equation*}\end{thm}

\beg{proof}   Let
$$h(t)= \sum_{i=1}^t \big(t\land \ell(t_i)-\ell(t_{i-1})\big)^+\si_{t_i}^{-1}J^\ell_{t_i} (J^\ell_T)^{-1} v,\ \ t\in [0,\ell(T)].$$ From  (H1) it is easy to see that $h\in \D(\ddd)$, where $(\ddd,\D(\ddd))$ is the    Malliavin divergence for the Brwonian motion $(W_t)_{t\in [0,\ell(T)]},$ see e.g. \cite{M,N}.
  Let $D_h$ be the Malliavin derivative along $h$.
Since $(V_t)_{t\ge 0}$ is independent of $(W_t)_{t\ge 0}$, we have $D_hV_t=0$, so that (\ref{2.10}) yields
\beq\label{D*} \d D_h X^\ell_t =(\nn_{D_h X^\ell_t}b_t)(X^\ell_t)\d t +\si_t \d h_{\ell(t)},\ \ D_h X^\ell_0=0.\end{equation} Then
$$D_h X^\ell_T= J^\ell_T \int_0^T(J^\ell_t)^{-1} \si_t \d h_{\ell(t)} = J^\ell_T \sum_{i=1}^n (J^\ell_{t_i})^{-1} \si_{t_i}\si_{t_i}^{-1} J^\ell_{t_i} (J^\ell_T)^{-1} \DD \ell(t_i) =\ell(T)v.$$
Therefore, \beq\beg{split} \label{2.1} P_T^\ell(\nn_v f) &=\E(\nn_vf)(X^\ell_T) =\ff 1 {\ell(T)} \E (\nn_{D_hX^\ell_T} f)(X^\ell_T) \\
&= \ff 1 {\ell(T)} \E\big\{D_h f(X^\ell_T)\big\} = \ff 1 {\ell(T)} \E\big\{f(X^\ell_T) \ddd(h)\big\}.\end{split}\end{equation}
To calculate $\ddd(h)$, let
$$h_{ik}(t)= (t\land \ell(t_i)-\ell(t_{i-1}))^+ e_k,\ \ F_{ik} =\big\<\si_{t_i}^{-1} J^\ell_{t_i} (J^\ell_T)^{-1} v, e_k\big\>$$ for $1\le i\le n, 1\le k\le d, t\in [0,\ell(T)],$ where $\{e_k\}_{k=1}^d$ is the canonical orthonormal basis on $\R^d$. Then
$$h(t)= \sum_{k=1}^d\sum_{i=1}^n F_{ik} h_{ik}(t),\ \ t\in [0,\ell(T)].$$
Noting that $h_{ik}$ is deterministic with $\int_0^{\ell(T)} |h'_{ik}(t)|^2\d t<\infty$, we have  $$\ddd(h_{ik})= \int_0^{\ell(T)} \<h_{ik}'(t),\d W_t\>= \<e_k, W_{\ell(t_i)}-W_{\ell(t_{i-1})}\>.$$   Thus, using the formula $\ddd(F_{ik}h_{ik})= F_{ik}\ddd(h_{ik})-D_{h_{ik}}F_{ik},$ we obtain
\beq\beg{split}\label{2.2} \ddd(h) &= \sum_{k=1}^d\sum_{i=1}^n \big\{F_{ik} \ddd(h_{ik}) -D_{h_{ik}}F_{ik}\big\}\\
&=\sum_{k=1}^d\sum_{i=1}^n \big\{F_{ik} \<e_k, W_{\ell(t_i)}-W_{\ell(t_{i-1})}\> -\<\si_{t_i}^{-1} D_{h_{ik}}(J^\ell_{t_i}(J^\ell_T)^{-1})v,e_k\>\big\}\\
&= \bigg\<\int_0^T\big(\si_t^{-1} J^\ell_t\big)^*\d W_{\ell(t)},  (J^\ell_T)^{-1} v\bigg\> - \sum_{k=1}^d\sum_{i=1}^n \big\<\si_{t_i}^{-1} D_{h_{ik}}(J^\ell_{t_i}(J^\ell_T)^{-1})v, e_k\big\>.\end{split}\end{equation} Since $\d h_{ik}(t)$ is supported on $(\ell(t_{i-1}), \ell(t_i))$ but $J^\ell_{t_i}$ is measurable with respect to
$\F_{\ell(t_{i-1})}:=\si\big\{W_{t}: t\le \ell(t_{i-1})\big\}$, we have $D_{h_{ik}} J^\ell_{t_i}=0$. So,
\beq\label{2.3} D_{h_{ik}} (J^\ell_{t_i}(J^\ell_T)^{-1}) = J^\ell_{t_i} D_{h_{ik}} (J^\ell_T)^{-1}= -J^\ell_{t_i} (J^\ell_T)^{-1} (D_{h_{ik}} J^\ell_T) (J^\ell_T)^{-1}.\end{equation} By (\ref{*1.2}), we have
$$\d D_{h_{ik}} J^\ell_t= (\nn_{D_{h_{ik}}X^\ell_t} \nn b_t)(X^\ell_t)J^\ell_t \d t + (\nn b_t)(X^\ell_t) D_{h_{ik}}J^\ell_t \d t,\ \ D_{h_{ik}}J^\ell_0=0.$$ Then
\beq\label{2.4}D_{h_{ik}}J^\ell_T = J^\ell_T \int_0^T (J^\ell_t)^{-1} (\nn_{D_{h_{ik}} X^\ell_t} \nn b_t)(X^\ell_t) J^\ell_t\d t.\end{equation} Moreover, by (\ref{D*}) we have
$$D_{h_{ik}} X^\ell_t= J^\ell_t \int_0^t (J^\ell_s)^{-1} \si_s \d h_{\ell(s)} = 1_{\{t_i\le t\}} (\ell(t_i)-\ell(t_{i-1})) J^\ell_t (J^\ell_{t_i})^{-1} \si_{t_i} e_k.$$
Combining this with (\ref{2.4}) we obtain
$$D_{h_{ik}} J^\ell_T= (\DD\ell(t_i))J^\ell_T \int_{t_i}^T (J^\ell_s)^{-1} \big(\nn_{J^\ell_s(J^\ell_{t_i})^{-1}\si_{t_i} e_k}\nn b_s\big)(X^\ell_s) J^\ell_s \d s.$$ Substituting this into (\ref{2.3}) we obtain
\beg{equation*}\beg{split} &\sum_{k=1}^d\sum_{i=1}^n \big\<\si_{t_i}^{-1} D_{h_{ik}} (J^\ell_{t_i}(J^\ell_T)^{-1})v,e_k\big\>\\
&= -\sum_{k=1}^d\sum_{i=1}^n (\DD\ell(t_i)) \si_{t_i}^{-1} J^\ell_{t_i}  \int_{t_i}^T \big\<(J^\ell_s)^{-1} \big(\nn_{J^\ell_s (J^\ell_{t_i})^{-1}\si_{t_i} e_k}\nn b_s\big)(X^\ell_s)J^\ell_s   (J^\ell_T)^{-1}v,e_k\big\>\d s \\
&=- \sum_{k=1}^d\int_0^T \d\ell(t) \int_t^T \big\<\si_t^{-1}J^\ell_t (J^\ell_s)^{-1} \big(\nn_{J^\ell_s (J^\ell_t)^{-1} \si_t e_k}\nn_{J^\ell_s(J^\ell_T)^{-1}v} b_s\big)(X^\ell_s), e_k\big\>\d s\\
&=-\int_0^T\d\ell(t)\int_t^T\sum_{k=1}^d \big\<\si_t^{-1}J^\ell_t (J^\ell_s)^{-1} \big(\nn\nn_{J^\ell_s(J^\ell_T)^{-1}v} b_s\big)(X^\ell_s)J^\ell_s (J^\ell_t)^{-1} \si_t e_k, e_k\big\>\d s\\
&=- \int_0^T \d \ell(t) \int_t^T {\rm Tr}\Big\{ \si_t^{-1} J_t^\ell(J_s^\ell)^{-1}\big(\nn\nn_{J_s^\ell(J_t^\ell)^{-1}v} b_s\big)(X^\ell_s)J_s^\ell(J_t^\ell)^{-1}\si_t \Big\}\d s.\end{split}\end{equation*}
Therefore, we derive from (\ref{2.2})  that $\ddd(h)= M_T^{\ell,v}$, and hence finish      the proof  by (\ref{2.1}).\end{proof}

\paragraph{Remark 2.1.} As indicated by Xicheng Zhang to the author that this result can also be proved using Remark 2.1 in \cite{WZ13}, which says that
$$P_T^\ell(\nn_vf)= \E\bigg(f(X_T^\ell) \sum_{i,k=1}^d\Big[D^*(h_k)\big\{(\nn X_T^\ell)^{-1}\big\}_{ki}-D_{h_k}\big\{(\nn X_T^\ell)^{-1}\big\}_{ki}\Big]v_i\bigg),$$ where $h_k\in\D(D^*)$ is such that
$$D_{h_k}X_T^\ell=\nn_{e_k} X_T^\ell.$$ Such $h_k$ can be constructed as that for  (3.8) in \cite{WXZ}. Noting that $\nn X_T^\ell=J_t^\ell,$ we can then prove Theorem \ref{T2.1} by calculating $D^*(h_k)$ and $D_{h_k}T_t^\ell$ as above.

\section{Proofs}
\beg{proof}[Proof of Theorem \ref{T1.1}] According to \cite[Theorem 2.4(1)]{W12}, the second assertion follows from the first. So, it suffices to prove the required integration by parts formula.  For any path $\ell$ of $S$ with $\ell(T)>0$, for any $\vv>0$, let
$$\ell_\vv(t)= \sum_{s\le t} \DD\ell(s) 1_{\{\DD\ell(s) \ge\vv\}},\ \ t\ge 0.$$ Then $\ell_\vv$ has finite many jumps on $[0,T].$ Moreover, $\d\ell_\vv(t)\to
\d\ell(t)$ on $[0,T]$ strongly as $\vv\to 0.$ Note that by (H1)
$$\big\|\si_t^{-1}J_t^{\ell_\vv}(J_T^{\ell_\vv})^{-1}\big\|+
 \int_t^T \big\|\si_t^{-1} J_t^{\ell_\vv} (J_s^{\ell_\vv})^{-1} \big(\nn_{J_s^{\ell_\vv}(J_t^{\ell_\vv})^{-1}\si_t e_k}\nn b_s\big)(X_s^{\ell_\vv}) J_s^{\ell_\vv}(J_t^{\ell_\vv})^{-1}\big\| \d s$$ is  bounded in $(t,\vv)\in [0,T]\times [0,1],$ and by \cite[Lemma 3.1]{WXZ}
$$\lim_{\vv\to 0} \E \sup_{t\in [0,T]} |X_t^{\ell_\vv} -X^\ell_t|^2 =0,$$ which together with (H1) implies
$$\lim_{\vv\to 0} \E \sup_{t\in [0,T]} \big(\|J_t^{\ell_\vv} -J_t^\ell\|^2+ \|(J_t^{\ell_\vv})^{-1}- (J_t^\ell)^{-1}\|^2\big)=0.$$ Combining these with (H1) and
(H2), we conclude that
\beg{equation*}\beg{split} &\lim_{\vv\to 0} P_T^{\ell_\vv} (\nn_vf) = \lim_{\vv\to 0} \E(\nn_vf)(X_T^{\ell_\vv}) = P_T^\ell (\nn_vf),\\
&\lim_{\vv\to 0} \E\big\{f(X_T^{\ell_\vv})M_T^{\ell_\vv,v}\big\}= \E\big\{f(X_T^{\ell})M_T^{\ell,v}\big\},\ \ f\in C_b^1(\R^d).\end{split}\end{equation*}
Therefore, first applying Theorem \ref{T2.1} to $\ell_\vv$ in place of $\ell$ then letting $\vv\downarrow 0$, we obtain
$$P_T^\ell(\nn_v f)=\ff 1 {\ell(T)} \E\big\{f(X_T^\ell)M_T^{\ell,v}\big\}$$ for all sample path $\ell$ of $S$ with $\ell(T)>0$. Since
$\E S(T)^{-\ff 1 2}<\infty$ implies $S(T)>0$, and noting that $X_T= X_T^S, M_T^v=M_T^{S,v}$, we obtain
\beq\label{D3} P_T^S(\nn_vf) = \ff 1 {S(T)} \E^S\big\{f(X_T)M_T^{v}\big\},\end{equation}where   $\E^S$ is the conditional expectation given $S$.
Moreover,  it follows from  (H1), (H2)  and $\E S(T)^{-\ff 1 2}<\infty$ that
\beq\label{D*2} \beg{split} &\E\Big|\ff {M_T^v} {S(T)} \Big|= \E\Big|\ff 1 {S(T)} \E^S M_T^v\Big| \\
&\le  \E \bigg[\ff 1 {S(T)}  \bigg(E^S \int_0^T |\si_t^{-1}J_t(J_T)^{-1}v|^2\d S(t)\bigg)^{1/2}\\
&\quad + \sum_{k=1}^d \ff 1 {S(T)}\int_0^T \d S(t) \int_t^T \|\si_t^{-1} J_t J_s^{-1}\|\cdot\big| \big(\nn_{J_s J_t^{-1}\si_t e_k}\nn_{J_sJ_t^{-1}v} b_s\big)(X_s)\big|\d s\bigg]\\
&\le |v| \Big(\ll_2(T)\e^{TK_1(T)}\E S(T)^{-\ff 1 2} + d T \ll_1(T)\ll_2(T) K_2(T) \e^{3TK_1(T)}\Big)<\infty.\end{split}\end{equation}   Then $  M_T^v\in L^1(S(T)^{-1}\d \P)$ and (\ref{D3}) yields that
$$P_T(\nn_vf)= \E P_T^S(\nn_vf) = \E\Big(\ff 1 { S(T)} f(X_T) M_T^v\Big).$$ This completes the proof.
\end{proof}

\beg{proof}[Proof of Corollary \ref{C1.2}] Assertion (1) follows immediately from (\ref{D*2}),   Theorem \ref{T1.1} and \cite[Theorem 2.4(1)]{W12} with $H(r)=r$.

Next,    by  (H1), (H2)  and the Burkholder inequality \cite[Theorem 2.3]{Z} (see also \cite[Lemma 2.1]{WXZ}),  for any $p>1$ there exists a constant $C(p)\ge 1$ such that
\beg{equation*}\beg{split}    \Big(\E  \ff {|M_T^v|^{\ff p{p-1}}}{ S(T)^{\ff p{p-1}}}\Big)^{\ff {p-1}p}
&\le   \bb(T)|v|+C(p) \Big(\E \ff {(\int_0^T|\si_t^{-1}J_t(J_T)^{-1}v|^2\d S(t))^{\ff p{2(p-1)}}}  {S(T)^{\ff p{p-1}}} \Big)^{\ff{p-1}p}\\
&\le \bb(T)|v|+ C(p)|v|   \ll_2(T)\e^{TK_1(T)}\big(\E S(T)^{-\ff p {2(p-1)}}\big)^{\ff {p-1}p}. \end{split}\end{equation*}    Then assertion (2) follows from \cite[Theorem 2.4(1)]{W12} with $H(r)= r^{\ff{p}{p-1}}$ and the fact that
$$|P_T(\nn_vf)|= \Big|\E\Big\{f(X^\ell_T)\ff{M_T^v}{S(T)}\Big\}\Big|\le (P_T|f|^p)^{\ff 1 p} \Big(\E  \ff {|M_T^v|^{\ff p{p-1}}}{ S(T) ^{\ff p{p-1}}}\Big)^{\ff {p-1}p},\ \ v\in\R^d.$$

Finally, by Theorem \ref{T1.1} and the Young inequality (see \cite[Lemma 2.4]{ATW09}), if $f\in C_b^1(\R^d)$ is non-negative then
\beq\label{2*} |P_T(\nn_vf)|= \Big|\E\Big\{\ff {f(X^\ell_T)M_T^v} {S(T)}\Big\}\Big| \le\dd {\rm Ent_{P_T}}(f) +\dd (P_T f) \log \E\exp\Big[\ff{M_T^v}{\dd S(T)}\Big],\ \dd>0.\end{equation} Obviously, by (H1) and (H2)
\beg{equation*}\beg{split} &\ff{M_T^v}{S(T)}\le  \bb(T)|v| +\ff 1 {S(T)} \int_0^T \<\si_t^{-1}J_tJ_T^{-1}v, \d W_{S(t)}\>,\\
& \E^S \exp\bigg[\ff 1 {\dd S(T)} \int_0^T \<\si_t^{-1}J_tJ_T^{-1}v, \d W_{S(t)}\>\bigg] \le \exp\bigg[\ff{\ll_2(T)^2|v|^2\e^{2TK_1(T)}}{2\dd^2S(T)}\bigg],\ \ \dd>0.\end{split}\end{equation*} Then
$$\log \E\exp\Big[\ff{M_T^v}{\dd S(T)}\Big]
 \le \ff{\bb(T)|v|}\dd + \log \E \exp\bigg[\ff{\ll_2(T)^2|v|^2\e^{2TK_1(T)}}{2\dd^2 S(T)} \bigg],\ \ \dd>0.$$ By combining this with (\ref{2*}) and \cite[Theorem 2.4(1)]{W12} for $H(r)= \e^{r/\dd}$, we prove (3).
\end{proof}

\beg{proof}[Proof of Corollary \ref{C1.3}] By \cite[Theorem 2.5(2)]{W12}, the second assertion follows from the first. So, we only need to prove the required shift-Harnack inequality \eqref{SH} for $v\ne 0$.  By Corollary \ref{C1.2}(3), we have
 $$|P_T(\nn_v f)|\le \dd {\rm Ent}_{P_T}(f) +   (P_T f)\Big(\bb(T)|v| + \dd \log \E\exp\Big[\ff{\ll_2(T)|v|^2}{2\dd^2 S(T)}\e^{\int_0^TK_1(t)\d t}\Big]\Big),\ \ \dd>0.$$
So, letting
$$\bb_v(\dd)= \bb(T)|v| + \dd \log \E\exp\Big[\ff{\ll_2(T)^2|v|^2\e^{2TK_1(T)}}{2\dd^2 S(T)}\Big],\ \ \dd>0,$$
we obtain from \cite[Proposition 2.3]{W12} that
\beq\label{*W} (P_Tf)^p\le (P_Tf^p(v+\cdot)) \exp\bigg[\int_0^1\ff{p}{1+(p-1)s} \bb_v\Big(\ff{p-1}{1+(p-1)s}\Big)\d s\bigg].\end{equation}
By the Jensen inequality,  for $\dd= \ff{p-1}{1+(p-1)s}$ we have
\beg{equation*}\beg{split} \E \exp\Big[\ff{\ll_2(T)^2|v|^2\e^{2TK_1(T)}}{2\dd^2 S(T)}\Big]&\le \bigg(\E\exp\Big[\ff{p^2\ll_2(T)^2|v|^2\e^{2TK_1(T)}}{2(p-1)^2S(T)}\Big]\bigg)^{\ff{(1+(p-1)s)^2}{p^2}}\\
&=\GG_{T,p}(|v|)^{\ff{(1+(p-1)s)^2}{p^2}}.\end{split}\end{equation*} Thus,
\beg{equation*}\beg{split} &\int_0^1\ff{p}{1+(p-1)s}\bb_v\Big(\ff{p-1}{1+(p-1)s}\Big)\d s\\
 &\le \bb(T)|v|\int_0^1\ff p{1+(p-1)s}\d s +\ff{p-1}p\log \GG_{T,p}(|v|)\d s\\
&= \ff{p\log p}{p-1} \bb(T)|v| + \ff{p-1}p\log\GG_{T,p}(|v|).\end{split}\end{equation*}
Then the proof is finished by combining this with (\ref{*W}).
\end{proof}

\beg{proof}[Proof of Corollary \ref{C1.4}] Since assertions in Corollaries \ref{C1.2} and \ref{C1.3} are unform in $V$, we may apply them for any deterministic path of $V$ in place of the process $V$, so that these two Corollaries remain true for $P_T^V$ in place of $P_T$, where
$$P_T^Vf(x) = \E^V (f(X_T(x)):= \E\big(f(X_T(x))|V\big).$$

Next, we observe that by the Markov property it suffices to prove the assertions for  $P_T^V$ in place of $P_T$ with $T\in (0,1]$.  In fact, for $T>1$ let
$$P_{1,T}^V f(x)= \E^V f(X_{1,T}(x)),\ \ f\in\B_b(\R^d), x\in\R^d,$$ where $(X_{1,t}(x))_{t\ge 1}$ solves the equation
$$X_{1,t}(x)= x + \int_1^tb_s(X_{1,s}(s))\d s +\int_1^t \si_s\d W_{S(s)}+V_t-V_1,\ \ t\ge 1.$$ Then by the Markov property of $X_t$ under $\E^V$, we obtain,
$$P_T^Vf= P_{1,T}^V(P_1^V f),\ \ f\in\B_b(\R^d).$$ Combining this with the assertions for $T=1$ and using the Jensen inequality, we prove the assertions for $T>1$. For instance, if for $p>1$ one has
$$|P_1^V(\nn f)|\le C(p) (P_1^V|f|^p)^{\ff 1 p},$$ then for any $T>1,$  \beg{equation*}\beg{split} &|P_T(\nn f)|= |\E P_{1,T}^VP_1^V(\nn f)|\le \E P_{1,T}^V|P_1^V(\nn f)|\\
&\le C(p) \E P_{1,T}^V (P_1^V|f|^p)^{\ff 1 p}\le C(p)(P_T|f|^p)^{\ff 1 p}= \ff{C(p)(P_T|f|^p)^{\ff 1 p}}{(1\land T)^{\ff  1 \aa}}.\end{split}\end{equation*} Below we prove assertions (1)-(3) for $T\in (0,1]$ respectively.

(1) Since $\bb(T)+ \ll_2(T)\e^{TK_1(T)}$ is bounded for $T\in (0,1]$, and by \cite[(ii) in the proof of Theorem 1.1]{WXZ}
$$\Big(\E S(T)^{-\ff{p}{2(p-1)}}\Big)^{\ff{p-1}p}\le \ff C{T^{\ff 1 \aa}},\ \ T\in (0,1]$$ holds for some constant $C>0$, the desired assertion follows from Corollary \ref{C1.2}(2).

(2) Let $\aa\in (1,2),$  and let $S_\aa$ be the subordinator induced by  the Bernstein function $r\mapsto r^{\ff\aa 2}$. Then as shown in \cite[Proof of Corollary 1.2]{WW} that
$$\E \ff 1 {S(T)^k}\le c_0\E\ff 1 {S_\aa(T)^k},\ \ k\ge 1, T\in (0,1]$$ holds for some constant $c_0\ge 1.$ Combining this with the third display from below  in the proof of \cite[Theorem 1.1]{GRW} for
$\kk=1$, i.e. (note the $\aa$ therein is $\aa/2$ here)
$$\E \e^{\ll/\tt S(t)} \le 1+\bigg(\exp\Big[\ff{c_1\ll^{\ff{\aa}{2(\aa-1)}}}{t^{\ff 1 {\aa-1}}}\Big]-1\bigg)^{\ff{2(\aa-1)}\aa}\le
\exp\Big[ \ff{c_2\ll}{t^{\ff 2 \aa}}+\ff{c_2\ll^{\ff\aa{2(\aa-1)}}}{t^{\ff 1{\aa-1}}}\Big],\ \ \ll,t\ge 0$$ for some constants $c_1,c_2>0,$ we obtain
\beq\label{*F}\beg{split}  \E\e^{\ll/S(T)} &\le 1+c_0\big(\E\e^{\ll/S_\aa(T)}-1\big)\le  \E\e^{c_0\ll/S_\aa(T)} \\
&\le \exp\Big[ \ff{c_3\ll}{t^{\ff 2 \aa}}+\ff{c_3\ll^{\ff\aa{2(\aa-1)}}}{t^{\ff 1{\aa-1}}}\Big],\ \ T\in (0,1], \ll\ge 0 \end{split}\end{equation} for some constant $c_3>0$.
By Corollary \ref{C1.2}(3) and (\ref{*F}), we prove the desired assertion.

(3) By (\ref{*F}), there exists a constant $c_4>0$ such that
$$\GG_{T,p}(r) \le \exp\Big[\ff{c_4 p^2r^2}{(p-1)^2T^{\ff 2\aa}}+\ff{c_4(pr)^{\ff\aa{\aa-1}}}{(p-1)^{\ff\aa{\aa-1}}T^{\ff 1{\aa-1}}}\Big],\ \ \ r\ge 0, T\in (0,1].$$ Then  there exists a constant $c_5>0$ such that
\beq\label{*F2} \beg{split}&\ff{p(\log p)\bb(T)|v|}{p-1}+ \ff{p-1}p \log \GG_{T,p}(|v|)\\
&\le \ff{c_5(p\log p)|v|}{p-1}
+ \ff{c_5 p|v|^2}{(p-1)T^{\ff 2\aa}}+\ff{c_5p^{\ff 1 {\aa-1}}|v|^{\ff\aa{\aa-1}}}{(p-1)^{\ff1{\aa-1}}T^{\ff 1{\aa-1}}},\ \ T\in (0,1], v\in\R^d.\end{split}\end{equation} By Corollary \ref{C1.3}, this implies the first inequality in (3) for some constant $C>0$.
Finally, the second inequality in (3) follows since (\ref{*F2}) and Corollary \ref{C1.3} imply
\beg{equation*}\beg{split} \sup_{x\in\R^d} \int_{\R^d} p_T(x,y)^{\ff p{p-1}}\d y &\le \bigg(\int_{\R^d} \exp\Big[-\ff{C(p\log p)|v|}{p-1}
-\ff{Cp|v|^2}{(p-1)T^{\ff 2\aa}}-\ff{Cp^{\ff 1{\aa-1}}|v|^{\ff\aa{\aa-1}}}{[(p-1)T]^{\ff 1{\aa-1}}}\Big]\d v\bigg)^{\ff {-1}{p-1}}\\
&\le \bigg(\int_{\{|v|\le T^{^\ff 1 \aa}\}} \d v\bigg)^{\ff{-1}{p-1}}\exp\Big[\ff{Cp(1+\log p)}{(p-1)^2}+\ff{Cp^{\ff 1 {\aa-1}}}{(p-1)^{\ff\aa{\aa-1}}}\Big]\\
&\le \ff 1 {T^{\ff d{\aa(p-1)}}}\exp\Big[\ff{C'p\log p}{(p-1)^2}+\ff{C'p^{\ff 1 {\aa-1}}}{(p-1)^{\ff\aa{\aa-1}}}\Big]\end{split}\end{equation*} for some constant $C'\ge C$.
\end{proof}

\paragraph{Acknowledgement.} The author would like to thank professor Xicheng Zhang for helpful comments.

\end{document}